\input amstex
\magnification=1200
\documentstyle{amsppt}
\NoBlackBoxes
\NoRunningHeads
\topmatter
\title Experimental detection of interactive phenomena and their analysis
\endtitle
\author Denis V. Juriev\endauthor
\affil ul.Miklukho-Maklaya 20-180 Moscow 117437 Russia\linebreak
(e-mail: denis\@juriev.msk.ru)\endaffil
\date math.GM/0003001\enddate
\abstract\nofrills The article is devoted to mathematical methods of
experimental detection of interactive phenomena in complex systems and 
their analysis.
\endabstract
\endtopmatter
\document
This article may be regarded as an account of conclusions from the ten year
author's practical researches on the joint of experimental mathematics,
experimental psychophysics and computer science in the problems of visual
perception in the interactive videosystems and from his parallel theoretical 
studies on interactivity and tactics. The current exposition is 
self-consistent and claims a minimal knowledge of the experimental or
theoretical basis, which underlies it. The goal of the article is to
discuss the methods of experimental detection of interactive phenomena and
of their analysis. Though the author had deal presumably with the interactive 
phenomena in perception the resulted scheme is applicable to many systems of
natural, behavioral, social and economical sciences and it is reasonable to
think that many concrete specialists will find it useful for their own needs 
in the nearest future.

\head I. Experimental detection of interactive phenomena\endhead

Let us consider a natural, behavioral, social or economical system $\Cal S$.
It will be described by a set $\{\varphi\}$ of quntities, which characterize
it at any moment of time $t$ (so that $\varphi=\varphi_t$). One may suppose
that the evolution of the system is described by a differential equation 
$$\dot\varphi=\Phi(\varphi)$$
and look for the explicit form of the function $\Phi$ from the experimental
data on the system $\Cal S$. However, the function $\Phi$ may depend on time,
it means that there are some hidden parameters, which control the system
$\Cal S$ and its evolution is of the form
$$\dot\varphi=\Phi(\varphi,u),$$
where $u$ are such parameters of unknown nature. One may suspect that such 
parameters are chosen in a way to minimize some goal function $K$, which may 
be an integrodifferential functional of $\varphi_t$:
$$K=K(\left[\varphi_{\tau}\right]_{\tau\le t})$$
(such integrodifferential dependence will be briefly notated as 
$K=K([\varphi])$ below). More generally, the parameters $u$ may be divided
on parts $u=(u_1,\ldots,u_n)$ and each part $u_i$ has its own goal function
$K_i$. However, this hypothesis may be confirmed by the experiment very 
rarely. In the most cases the choice of parameters $u$ will seem accidental
or even random. Nevertheless, one may suspect that the controls $u_i$ are 
{\sl interactive}, it means that they are the couplings of the pure controls 
$u_i^\circ$ with the {\sl unknown or incompletely known\/} feedbacks:
$$u_i=u_i(u_i^\circ,[\varphi])$$
and each pure control has its own goal function $K_i$. Thus, it is
suspected that the system $\Cal S$ realizes an {\sl interactive game}.
There are several ways to define the pure controls $u_i^\circ$. One of them
is the integrodifferential filtration of the controls $u_i$:
$$u^\circ_i=F_i([u_i],[\varphi]).$$
To verify the formulated hypothesis and to find the explicit form of the
convenient filtrations $F_i$ and goal functions $K_i$ one should use the
theory of interactive games, which supplies us by the predictions of the
game, and compare the predictions with the real history of the game for
any considered $F_i$ and $K_i$ and choose such filtrations and goal functions,
which describe the reality better. One may suspect that the dependence of
$u_i$ on $\varphi$ is purely differential for simplicity or to introduce the
so-called {\sl intention fields}, which allow to consider any interactive
game as differential. Moreover, one may suppose that
$$u_i=u_i(u_i^\circ,\varphi)$$
and apply the elaborated procedures of {\sl a posteriori\/} analysis and
predictions to the system.

In many cases this simple algorithm effectively unravels the hidden 
interactivity of a complex system.

\head II. Analysis of interactive phenomena\endhead

Below we shall consider the complex systems $\Cal S$, which have been yet
represented as the $n$-person interactive games by the procedure described
above. 

\subhead 2.1. Functional analysis of interactive phenomena\endsubhead
To perform an analysis of the interactive control let us note that often for 
the $n$-person interactive game the interactive controls 
$u_i=u_i(u_i^\circ,[\varphi])$ may be represented in the form 
$$u_i=u_i(u_i^\circ,[\varphi];\varepsilon_i),$$
where the dependence of the interactive controls on the arguments
$u_i^\circ$, $[\varphi]$ and $\varepsilon_i$ is known but the 
$\varepsilon$-parameters $\varepsilon_i$ are the unknown or incompletely
known functions of $u_i^\circ$, $[\varepsilon]$. Such representation is
very useful in the theory of interactive games and is called the 
{\sl $\varepsilon$-representation}. 

One may regard $\varepsilon$-parameters as new magnitudes, which characterize
the system, and apply the algorithm of the unraveling of interactivity to
them. Note that $\varepsilon$-parameters are of an existential nature 
depending as on the states $\varphi$ of the system $\Cal S$ as on the
controls. 

The $\varepsilon$-parameters are useful for the functional analysis of
the interactive controls described below.

First of all, let us consider new integrodifferential filtrations $V_\alpha$:
$$v^\circ_\alpha=V_\alpha([\varepsilon],[\varphi]),$$
where $\varepsilon=(\varepsilon_1,\ldots,\varepsilon_n)$. 
Second, we shall suppose that the $\varepsilon$-parameters are expressed via 
the new controls $v^\circ_\alpha$, which will be called {\it desires:}
$$\varepsilon_i=\varepsilon(v^\circ_1,\ldots,v^\circ_m,[\varphi])$$
and the least have the goal functions $L_\alpha$. The procedure of unraveling
of interactivity specifies as the filtrations $V_\alpha$ as the goal functions 
$L_\alpha$.

\remark{Example} Let us considered the interactive videosystem directed by
the eye movements of an observer. The pure controls are the slow movements
of eyes, whereas saccads are considered as a result of the unknown feedbacks
(tremor is supposed to be random). Many classical and modern experiments
clarifies the role of saccads in the formation of the stable and complete
final image so such formation may be regarded as their goal function. The 
functional analysis of the eye movements extracts the parameters (the normal 
forms), which describe saccads in the concrete interactive videosystems. The 
normal forms are extremely interesting in the multi-user mode when the saccads 
of various observers begin to be correlated and synchronized.
\endremark

\subhead 2.2. The second quantization of desires\endsubhead
Intuitively it is reasonable to consider systems with a variable number
of desires. It can be done via the second quantization. 

To perform the second quantization of desires let us mention that they
are defined as the integrodifferential functionals of $\varphi$ and
$\varepsilon$ via the integrodifferential filtrations. So one is able
to define the linear space $H$ of all filtrations (regarded as classical 
fields) and a submanifold $M$ of the dual $H^*$ so that $H$ is naturally
identified with a subspace of the linear space $\Cal O(M)$ of smooth functions
on $M$. The quantized fields of desires are certain operators in the
space $\Cal O(M)$ (one is able to regard them as unbounded operators in its
certain Hilbert completion). The creation/annihilation operators are
constructed from the operators of multiplication on an element of $H\subset
\Cal O(M)$ and their conjugates.

To define the quantum dynamics one should separate the quick and slow time.
Quick time is used to make a filtration and the dynamics is realized in
slow time. Such dynamics may have a Hamiltonian form being governed by
a quantum Hamiltonian, which is usually differential operator in $\Cal O(M)$.

If $M$ coincides with the whole $H^*$ then the quadratic part of a Hamiltonian
describes a propagator of the quantum desire whereas the highest terms
correspond to the vertex structure of self-interaction of the quantum field. 
If the submanifold $M$ is nonlinear the extraction of propagators and 
interaction vertices is not straightforward.

\subhead 2.3. SD-transform and SD-pairs\endsubhead
The interesting feature of the proposed description (which will be called the
{\it S-picture}\/) of an interactive system $\Cal S$ is that it contains as 
the real (usually personal) subjects with the pure controls $u_i$ as the 
impersonal desires $v_\alpha$. The least are interpreted as certain 
perturbations of the first so the subjects act in the system by the 
interactive controls $u_i$ whereas the desires are hidden in their actions. 

One is able to construct the dual picture (the {\sl D-picture\/}),
where the desires act in the system $\Cal S$ interactively and the
pure controls of the real subjects are hidden in their actions.
Precisely, the evolution of the system is governed by the equations
$$\dot\varphi=\tilde\Phi(\varphi,v),$$
where $v=(v_1,\ldots,v_m)$ are the $\varepsilon$-represented interactive 
desires:
$$v_\alpha=v_\alpha(v^\circ_\alpha,[\varphi];\tilde\varepsilon_\alpha)$$
and the $\varepsilon$-parameters $\tilde\varepsilon$ are the unknown or
incompletely known functions of the states $[\varphi]$ and the pure
controls $u_i^\circ$.

D-picture is convenient for a description of systems $\Cal S$ with a
variable number of acting persons. Addition of a new person does not
make any influence on the evolution equations, a subsidiary term to
the $\varepsilon$-parameters should be added only.

The transition from the S-picture to the D-picture is called the
{\it SD-transform}. The {\it SD-pair\/} is defined by the evolution
equations in the system $\Cal S$ of the form
$$\dot\varphi=\Phi(\varphi,u)=\tilde\Phi(\varphi,v),$$
where $u=(u_1,\ldots,u_n)$, $v=(v_1,\ldots,v_m)$, 
$$\aligned
u_i=&u_i(u_i^\circ,[\varphi];\varepsilon_i)\\
v_\alpha=&v_\alpha(v^\circ_\alpha,[\varphi];\tilde\varepsilon_\alpha)
\endaligned$$
and the $\varepsilon$-parameters $\varepsilon=(\varepsilon_1,\ldots,
\varepsilon_n)$ and $\tilde\varepsilon=(\tilde\varepsilon_1,\ldots,
\tilde\varepsilon_m)$ are the unknown or incompletely known functions of
$[\varphi]$ and $v^\circ=(v^\circ_1,\ldots,v^\circ_m)$ or
$u^\circ=(u^\circ_1,\ldots,u^\circ_n)$, respectively. 

Note that the S-picture and the D-picture may be regarded as complementary
in the N.Bohr sense. Both descriptions of the system $\Cal S$ can not be 
applied to it simultaneously during its analysis, however, they are compatible 
and the structure of SD-pair is a manifestation of their compatibility.
The choice of a picture is an action of our {\it attention:\/} it is 
concentrated on the personal subjects in S-picture {\it (the self-conscious 
attention)\/} whereas it is concentrated on the impersonal desires in 
D-picture {\it (the creative attention)}.

\subhead 2.4. Verbalization of SD-pairs and synlinguism\endsubhead
The main problem is to interrelate the S- and D-pictures of the system 
$\Cal S$. One way is a {\it verbalization\/} of SD-pairs. Let us remind
a definition of the verbalizable interactive game.
 
An interactive game of the form
$$\dot\varphi=\Phi(\varphi,u)$$
with $\varepsilon$--represented couplings of feedbacks 
$$u_i=u_i(u^\circ_i,[\varphi];\varepsilon_i)$$
is called {\it verbalizable\/} if there exist {\sl a posteriori\/}
partition $t_0\!<\!t_1\!<\!t_2\!<\!\ldots\!<\!t_n\!<\!\ldots$ and the 
integrodifferential functionals
$$\aligned
\omega_n&(\vec\varepsilon(\tau),\varphi(\tau)|
t_{n-1}\!\leqslant\!\tau\!\leqslant\!t_n),\\
u^*_n&(u^\circ(\tau),\varphi(\tau)|
t_{n-1}\!\leqslant\!\tau\!\leqslant\!t_n)
\endaligned$$ 
such that
$$\omega_n=\Omega(\omega_{n-1},u^*_n;\varphi(\tau)|
t_{n-1}\!\leqslant\!\tau\!\leqslant\!t_n),$$ 
quantities $\omega_n$ are called the {\it words}. 

Let us now consider the SD-pair and suppose that both S- and
D-pictures are verbalizable with the {\sl same\/} $\omega_n$.
The fact that $\omega_n$ are the same for both S- and D-pictures
is called their {\it synlinguism}. One may characterize it poetically by the
phrase that {\sl ``the speech of real subjects is resulted in the same text 
as a whisper of the impersonal desires''}. The existential character of the 
synlinguism should be stressed. Really it is not derived from the fact that 
the objective states $\varphi$ of the system $\Cal S$ are the same in the 
S- and D-pictures. The synlinguism interrelates the different 
$\varepsilon$-parameters of existential nature in both pictures.

The synlinguism is very important in the analysis of tactical phenomena, 
which essentially used the concept of verbalization in their definition.
To my mind the synlinguism lies in the basis of psychophysical nature of
mutual understanding of the independent subjects of a dialogue communication.
In this situation it allows to identify the personal interpretations with
the impersonal ones, unraveling the role of impersonal desires as bearers 
of the objective sense and its dynamics.

To the verbalizable SD-pairs some procedures of linguistic analysis are 
applicable. Some of them are inherited from the verbalizable interactive 
games (the grammatical analysis), some are specific (the explication and 
analysis of objective sense).

\head III. Conclusions\endhead

Thus, mathematical procedures of the experimental detection of interactive
phenomena in complex natural, behavioral, social and economical systems and
their analysis are described. The special attention is concentrated on the 
role of desires and their second quantization as well as on the abstract 
structure of SD-pairs, their verbalization and the synlinguism.
\enddocument